\documentclass{elsart}
\usepackage{amssymb}
\newtheorem{remark}{Remark}[section]
\newtheorem{definition}{Definition}[section]
\newtheorem{lemma}{Lemma}[section]
\newtheorem{theorem}{Theorem}[section]

\begin{document}
\begin{frontmatter}

\title{Generalized momenta of mass and their applications to the flow of compressible fluid}
\author{Olga Rozanova }
\thanks{Supported by DFG Project 436 RUS 113/823/0-1}
\address { Mech.\&Math. Faculty, Moscow State University, Moscow,
119992, Russia }

\ead {rozanova@mech.math.msu.su}

\begin{abstract}We present a technique that allows to obtain certain results in the compressible
fluid theory: in particular, it is a nonexistence result for the
highly decreasing at infinity solutions to the Navier-Stokes
equations, the construction of the solutions with uniform
deformation and the study of behavior of the boundary of a material
volume of liquid.
\end{abstract}

\begin{keyword}
compressible fluid \sep the Cauchy problem \sep exact solutions \sep
material volume 
\MSC 53Q \sep 53L65
\end{keyword}
\end{frontmatter}

%







%
%

%
%
\section{Main equations}
The motion of compressible viscous, heat-conductive, 
polytropic fluid in ${\mathbb R}\times{\mathbb R}^n,\, n\ge 1,$ is
governed by the compressible Navier-Stokes (NS) equations
$$\partial_t \rho+{\rm div}_x (\rho v)=0,\eqno(1.1)$$
$$\partial_t(\rho v)+{\rm div}_x (\rho v \otimes v)+\nabla_x p={\rm Div} T,\eqno(1.2)$$
$$\partial_t\left(\frac{1}{2}\rho |v|^2+\rho e\right)+{\rm div}_x \left((\frac{1}{2}\rho |v|^2+\rho e+p)v
\right)={\rm div}(Tv)+k\Delta_x \theta, \eqno(1.3)     $$ where
$\rho, \, u=(u_1,...,u_n),\, p, \, e,\,\theta \,$ denote the
density, velocity, pressure, internal energy and absolute
temperature, respectively; $\,T\,$ is the stress tensor given by the
Newton law $T=T_{ij}=\mu \,(\partial_iv_j+\partial_j v_i)+\lambda
\,{\rm div} u\, \delta_{ij},$  with constant coefficients of
viscosity $\mu$ and $\lambda$
 ($\mu\ge 0, \,\lambda+\frac{2}{n}\mu\ge
0$), $k\ge 0$ is the coefficient of heat conduction. We denote ${\rm
Div}$ and $\rm div$ the divergency of tensor and vector,
respectively. The state equations are $\,p=(\gamma-1)\rho e,\,$ and
$p=R\rho \theta,$ where $\gamma>1$ and $R>0$ are  the specific heat
ratio and the universal gas constant, respectively. Thus, we can
consider (NS) as a system for unknown $\rho,\,u,\,p.$ Indeed, (NS)
and the state equations give
$$\partial_t p +
(v,\nabla_x p)+\gamma p\,\,{\rm
div}v=(\gamma-1)\,\sum\limits_{i,j=1}^n\,T_{ij}
\partial_j v_i+\frac{k}{R}\,\Delta\frac{p}{\rho}.\eqno(1.4)$$
For classical solutions (NS) is equivalent to system (1.1, 1.2,
1.4), denoted (NS*) for short.

(NS*) is supplemented with the initial data
$$(\rho,v,p)\Bigl|_{t=0}=(\rho_0(x),v_0(x),p_0(x))\in C^2({\mathbb
R}^n).$$ If $\mu=\lambda=k=0,$ we get the gas dynamic (GD)
equations.

\section{Conservation laws, generalized momenta of mass and the decay rate}

System (NS) is the differential form of conservation laws for the
material volume ${\mathcal V}(t)$ (i.e. the volume, consisting from
the same particles). (NS) expresses  conservation of mass
$m=\int\limits_{{\mathcal V}(t)}\rho \, d x,$ momentum
$P=\int\limits_{{\mathcal V}(t)}\rho v\, d x,$ and  total energy
$\mathcal E=\int\limits_{{\mathcal V}(t)}\left(\frac{1}{2}\rho
|v|^2+\rho e\right)\, d x \,=E_{k}(t)+E_{i}(t).$ Here $E_k(t)$ and
$E_i(t)$ are the kinetic and internal components of energy,
respectively.

When ${\mathcal V}(t)={\mathbb R}^n,$ the conservation of mass,
angular momentum and energy take place provided the components of
solution  decrease at infinity sufficiently quickly.

\begin{definition}
We say that the classical solution solution to (NS) belongs to the
class $K(M(t),\alpha),$ if
 there exist a positive vector-function $M(t)=(M_v(t),M_\rho(t),$ $M_p(t))$ with components from
$C([0,\infty)),$  a constant vector $\alpha=(\alpha_v,\alpha_{Dv},
\alpha_\rho,\alpha_p,\alpha_\theta)$ and  constants $R_0>0,\,T\ge
0,$ such that $|v(t,x(t))|\le M_v(t) |x(t)|^{\alpha_v},\,
|Dv(t,x(t))|\le M_{Dv}(t) |x(t)|^{\alpha_{Dv}},\,$ $\rho(t,x(t))\le
M_\rho(t) |x(t)|^{\alpha_\rho},$ $\,p(t,x(t))\le M_p(t)
|x(t)|^{\alpha_p},\,\theta(t,x(t))\le M_{\theta}(t)
|x(t)|^{\alpha_{\theta}}\,$ for all trajectories $x(t)$ such that
$|x(t)|>R_0,\, t>T.$
\end{definition}

Let us introduce a functional
$$G_\phi(t)=\int\limits_{{\mathcal
V}(t)}\rho(t,x)\phi(|{x}|)\,dx.$$ If $\phi(|{ x}|)=\frac{1}{2} |{
x}|^2,$ then $G_\phi(t)$ is the usual momentum of mass. By analogy
for others $\phi $ we call $G_\phi(t)$ the generalized momentum of
mass.

\subsection{Decay rate for (NS) equations}

Let us choose the appropriate class $K(M(t),\alpha)$ to guarantee
conservation of mass, momentum, energy and convergence of the mass
momentum $G(t)$ for (NS) system. It is sufficient to set $\alpha
=(\alpha_v, \alpha_{Dv},\alpha_\rho,
\alpha_p,\alpha_\theta)=(-n,-n-1,-n-2-\varepsilon,-n-\varepsilon,-n)$
with a constant $\varepsilon>0.$ We denote this class $K_{NS}.$

If the heat conductivity is zero, we do not need to require the
decay of $\theta,$ i.e
$\alpha=(-n,-n-1,-n-2-\varepsilon,-n-\varepsilon,\alpha_\theta)$
with an arbitrary last component. We denote this class $K_{NS_0}.$

\subsection{Decay rate for (GD) equations}

In this case the behavior of velocity  is less restrictive. Here it
is sufficient to set $\alpha =(\alpha_v,
\alpha_{Dv},-n-2-\varepsilon,-n-\varepsilon,\alpha_\theta)$ with
$\varepsilon>0,\, \alpha_v\le 1$ and arbitrary $\alpha_{Dv}$ and
$\alpha_\theta.$    The components of velocity may rise as $|{\bf
x}|\to \infty$. We  denote this class $K_{GD}.$

\section{Behavior of generalized momenta of mass on solutions}

We denote $\,{\bf \sigma}=(\sigma_1,...,\sigma_K)$ the vector with
components $\,\sigma_k=v_i x_j-v_j x_i,$
$\,i>j,\,i,j=1,...,n,\,k=1,...,K,\,K={\rm C}_n^2.$

The following lemma describes the behavior of generalized momenta
for (GD) system.

\begin{lemma} Let us suppose that  $\phi (|{x}|)$
belongs to the class $C^2$ inside ${\mathcal V}(t)$. For classical
solution to system (GD) the following equalities take place:
$$G'_\phi(t)=\int\limits_{{\mathcal V}(t)}
\frac{\phi'(|{x}|)}{|{ x}|} ({ v},{ x})\rho\,dx,$$
$$G''_\phi(t)=I_{1,\phi}(t)+I_{2,\phi}(t)+I_{3,\phi}(t)+I_{4,\phi}(t),$$
where
$$I_{1,\phi}(t)=\int\limits_{{\mathcal V}(t)}\frac{\phi''(|{ x}|)}{|{
x}|^2} |({ v},{ x})|^2\rho\,dx,\qquad
I_{2,\phi}(t)=\int\limits_{{\mathcal V}(t)}\frac{\phi'(|{ x}|)}{|{
x}|^3} |{ \sigma}|^2\rho\,dx,$$
$$I_{3,\phi}(t)=\int\limits_{{\mathcal V}(t)}(\phi''(|{ x}|)+
(n-1)\frac{\phi'(|{ x}|)}{|{ x}|} ) p\,dx,\quad
I_{4,\phi}(t)=-\int\limits_{\partial{\mathcal V}(t)} \frac{\phi'(|{
x}|)}{|{ x}|} ({ x},{\nu}) p\,d\partial{\mathcal V},$$ where $\nu$
is a unit outer normal to $\mathcal V.$
\end{lemma}

The {\it proof} is a direct application of the Stokes formula.

Let us consider $ \phi(|{ x}|)=\frac{1}{2}|{ x}|^2$ and $ {\mathcal
V}(t)={\mathbb R}^n.$ We denote the respective functional
$G_\phi(t)$ as $G(t).$ In this case the following lemma is true.

\begin{lemma}
For  classical solutions to (NS),  $k=0,$  of the class $K_{NS_0}$
  and solutions to (GD) of the class $K_{GD}$
the following equality holds:
$$G''(t)= 2E_k(t)+n(\gamma-1)E_i(t).\eqno(3.1)
$$
\end{lemma}
\begin{pf}
First of all, in (NS) case $\,G''(t)=
I_1(t)+I_{2}(t)+I_{3}(t)+I_4(t)+I_5(t),$ where
$I_{1}(t)=\int\limits_{{\mathbb R}^n}\frac{|({ u},{ x})|^2 }{|{
x}|^2} \rho\,dx,\, I_{2}(t)=\int\limits_{{\mathbb R}^n}\frac{|{\bf
\sigma}|^2} {|{ x}|^2} \rho\,dx,\, I_{3}(t)=n\int\limits_{{\mathbb
R}^n}p\,dx,\, I_{4,\phi}(t)=
-\lim\limits_{R\to\infty}\int\limits_{\partial{ B_R}}  ({ x},{\nu})
p\,d{\partial B_R},\,I_5(t)=\int\limits_{{\mathbb R}^n}({\rm Div}
T,x)\,dx\,=\,\lim\limits_{R\to\infty}\int\limits_{\partial
B_R}T_{ij}x_i\nu_j \, d \partial B_R-(2\mu+n\lambda)\int\limits_{
B_R}{{\rm div} v}\,dx\,
=\,\lim\limits_{R\to\infty}\,\int\limits_{\partial B_R}(T_{ij}x_i-
(2\mu+n\lambda)v_i\,\delta_{ij}) \nu_j \,d
\partial B_R,\,$ where $B_R:|x|\le R,$ we sum the repeated
indices. Owing to respective assumptions on the decay of solution as
$|{\bf x}|\to \infty $ all improper integrals converge and
$I_4(t)=I_5(t)=0.$ To end the proof we note that
$\,I_{1}(t)+I_{2}(t)=2E_k(t),\,I_{3}(t)=n(\gamma-1)E_i(t).$ In (GD)
case the integral $I_5$ is missing.
\end{pf}

\begin{remark} Equality (3.1) for (GD) was obtained firstly in \cite{Chemin}.
\end{remark}
\begin{remark} It follows from (3.1) that $n(\gamma-1){\mathcal E}\le G''(t)\le 2{\mathcal E}. $
The inequality implies the same two sided estimate of $G(t)$ that
was obtained in 3D for the sum of squared semi-axis  of a
non-rotating gas ellipsoid next to vacuum \cite{PMM}.
\end{remark}

\section{Necessary conditions for existence of global solution}

If the initial density $\rho_0$ is compact, then in arbitrary space
dimensions no solution to (NS) from $C^1([0,\infty), H^m({\mathbb
R}^n)),\,m\ge\left[\frac{n}{2}\right]+2,$ exists \cite{Xin}. This
blowup result depends crucially on the assumption about  compactness
of support of the initial density. Thus, the question remains: is it
true that the global in $t$ smooth solution exists for any smooth
initial data in the case where the support of initial density
coincides with the whole space? Below we find necessary conditions
of existence of the global solution with prescribed decay rate  as
$|{x}|\to \infty$.

\begin{theorem}
If the global in $t$ solution to (NS) of the class $K_{NS_0}$
exists, then the solution components grow as $t\to\infty $ at least
such that
$$\int_0^t\,M_v(\tau)\,d\tau=O(t^{1-\alpha_v}),\quad
M_\rho(t)=O(t^{2+\epsilon}),\quad \epsilon>0,\quad t\to\infty.
\eqno(4.1)$$
\end{theorem}

\begin{pf}
It follows from Lemma 2 that $G''(t)\ge (\gamma-1){\mathcal E},$
$$
G(t)\ge \frac{(\gamma-1)n{\mathcal E}}{2}t^2+G'(0)
t+G(0).\eqno(4.2)$$ Let us get the  estimate of $G(t)$ from above.

 We consider a material volume $\mathcal V(t)$ that initially
coincides with the ball $|{\bf x}|\le R_0$ (see Definition 1). We
denote $B_{R(t)}=(x:|{ x}|\le R(t))$ a ball that contains $\mathcal
V(t).$ From the definition of  $K_{NS}$ we have
$$
\frac{d R(t)}{dt}\le R^{\alpha_v}(t)\, M_v(t),\quad R(t)\le
\left((1-\alpha_v) \int_0^t\,M_v(\tau)\,d\tau
+R_0^{1-\alpha_v}\right)^{\frac{1}{1-\alpha_v}}.
$$
Further, $$\,G(t)=\frac{1}{2}\,\left(\int\limits_{B_{R(t)}} \rho |{
x}|^2 \,dx+ \int\limits_{{\mathbb R}^n \backslash B_{R(t)}} \rho |{
x}|^2 \,dx\right)\le \frac{R^2(t)}{2}\,m+$$
 $$\frac{1}{2}\,M_\rho(t)\,\int\limits_{{\mathbb R}^n \backslash
B_{R(t)}}  |{ x}|^{2+\alpha_\rho} \,dx\le$$
$$
O\left(\left(\int_0^t\,M_v(\tau)\,d\tau
\right)^{\frac{2}{1-\alpha_v}}
\right)+O\left(M_\rho(t)\,\left(\int_0^t\,M_v(\tau)\,d\tau
\right)^{\frac{-\epsilon}{1-\alpha_v}} \right). $$ If the growth of
$M_v(t)$ and $M_\rho(t)$ is less then  prescribed in (4.1), then the
latter inequality contradicts to (4.2).
\end{pf}

\begin{remark} In particulary, Theorem 1 implies that there exist
no smooth solution to (NS) from the class $K_{NS_0}$ with
$M_v(t)=\rm const$ and/or $M_\rho(t)=\rm const.$
\end{remark}

The analog of Theorem 1 in the case of (GD) equations is the
following.

\begin{theorem} If the global in $t$ solution to (GD) from the class $K_{GD}$
exists  then the functions $M_\rho(t)$ are $M_v(t),$ restricting the
density and velocity behavior  as $ t\to\infty,$  grow at least such
that
 $M_\rho(t)=O(t^{2+\epsilon}),\,
\epsilon>0,\,$ and $\,\int_0^t\,M_v(\tau)\,d\tau=O(\ln t)$ (for
$\alpha_v=1$) or as prescribed in (4.6) (for $\alpha_v<1$).
\end{theorem}

The statement can be proved exactly as Theorem 1.


\section{Motion with uniform deformation}
We dwell in this section on (GD) equations, however due to the
specific choice of the velocity field the results remain true for
(NS) with $k=0.$

 The motion with uniform
deformation (i.e with linear profile of velocity $v(t,x)=A(t)\,{ x},
\,A(t)$ is a matrix ($n\times n$)) was considered in many works. Let
us mention \cite{Ovs56}, where this special solution to (GD) was
firstly constructed in 3D and \cite{Bogoyavlenski}, where many
applications are given. Generally speaking, for this solution the
velocity, density and pressure may blow up as $t\to T<\infty$ and
$|{\bf x}|\to \infty.$ Below we show how the solution with uniform
deformation can be constructed by means of the generalized momenta
of mass. Moreover, the requirement of finiteness of energy and mass
momentum prohibits the blow up, therefore the respective solution is
globally in time smooth.

Below we compare the procedure of constructing of the simplest
solution with the linear velocity profile
$$v(t,x)=a(t)\,{x}\eqno(5.1)$$
for a different choice of generalized momenta.

\subsection{Basing on the momentum of mass}
Let us consider once more the particular case $\phi(|{
x}|)=\frac{1}{2}|{ x}|^2,$ the respective functional is the usual
momentum of mass, $G(t)$.

Lemma 2 and (1.4) imply
$$G'(t)=2 a(t) G(t),\qquad  a'(t)=-a^2(t)+K G^{-\frac{(\gamma-1)n+2}{2}} (t),  \eqno(5.2)$$
with a constant $K>0$  depending on initial data.
 If the initial data satisfy the compatibility condition
 $$ p_0'(|{ x}|)= - (\gamma-1)G^{-1}(0)E_i(0)\rho_0(|{ x}|),\eqno(5.3)$$
then the density and pressure can be found from (1.1) and (1.4) as
$$\rho(t,|x|)=\exp(-n\int_0^t a(\tau) d\tau)
\rho_0(|x|\exp(-\int_0^t a(\tau) d\tau)),\eqno(5.4)$$
$$p(t,|x|)=\exp(-n\gamma\int_0^t a(\tau) d\tau)
\rho_0(|x|\exp(-\int_0^t a(\tau) d\tau)).\eqno(5.5)$$ It follows
from (5.2) that $a(t)$ solves
$$a'(t)=-a^2(t)+ K_1
\exp\left(-((\gamma-1)n+2)\int_0^t\,a(\tau)\,
d\tau\right),\eqno(5.6)$$ with $K_1=K(G(0))^{-((\gamma-1)n+2)/2}.$
Since $K_1>0,$ the solution $a(t)$ remains bounded for all $t\ge 0$
and $a(t)=O(t^{-1}),\,t\to\infty$ \,\cite{RozHyp}.

\begin{remark}
 In the case of globally
smooth solution discussed above, the velocity  $v(t,x)=a(t)\,{ x},$
therefore $\alpha_v=1.$ Here $M_v(t)=a(t)=O(t^{-1}),$ and
$\int_0^t\,M_v(\tau)\,d\tau=O(\ln
t),\,M_\rho(t)=O(t^{2+\epsilon}),\, t\to\infty.$ Thus, according
Theorem 1, the solution presents a less possible rate of its
component in $t.$

\end{remark}
\begin{remark} In \cite{RozNova} the solutions with linear profile of velocity
for (GD) equations (including the presence of the Coriolis force and
damping) were constructed for a general matrix $A(t)$ in 2D and for
particular cases in 3D.
\end{remark}

\subsection{The "excluding pressure" case}
Let us choose as $\phi$ a function, proportional to the fundamental
solution of the Laplace operator, namely
$$\phi(|{x}|)=\ln |{x}|,\quad n=2, \qquad \phi(|{ x}|)=
|{ x}|^{2-n},\quad n\ne 2.
$$
In this case Lemma 1 implies
$$G'_\phi(t)=\lambda_1 (n) \int\limits_{{\mathbb R}^n}
\frac{({v},{ x})}{|{ x}|^n} \rho\,dx,$$
$$G''_\phi(t)=\lambda_2(n)\,\int\limits_{{\mathbb R}^n}\frac{|({ v},{ x})|^2}{|{ x}|^{n+2}}\rho\,dx+
\lambda_3(n)\,\int\limits_{{\mathbb R}^n}\frac{|{\bf
\sigma}|^2\rho}{|{ x}|^{n+2}}\,dx +
\frac{p(t,0)}{\omega_{n-1}(2-n)},$$ where
$$\lambda_1(n)= 1, \quad  n=1,2 \quad \mbox{and}\quad
2-n, \quad n\ge 3,$$
$$\lambda_2(n)= 0, \quad  n=1; \quad -1, \quad n=2, \quad \mbox{and}\quad
(1-n)(2-n), \quad n\ge 3,$$
$$\lambda_3(n)= 0, \quad  n=1; \quad 1, \quad n=2, \quad \mbox{and}\quad
2-n, \quad n\ge 3.$$ Here we take as ${\mathcal V}(t)$ the space
${\mathbb R}^n, $ all improper integrals are supposed convergent
both at the origin and at infinity. Actually, this signifies that
$v(t,0)=0.$

We use the value of pressure only in the origin; in this sense the
pressure in the remaining space is excluded.

Let us consider the case $\, n\ge 3.$ Here
 $\,G'_{\phi}(t)=\lambda_1(n) \,a(t)G_{\phi}(t),$$$ G''_{\phi}(t)=\lambda_2(n)
\,a^2(t)\,G_{\phi}(t)+
\frac{p(0,0)}{\omega_{n-1}(2-n)}\,\exp\left(-\gamma
n\,\int_0^t\,a(\tau)\,d\tau\right).$$ The latter system implies
$$a'(t)=-a^2(t)+ K_2 \exp\left(-((\gamma-1)n+2)\int_0^t\,a(\tau)\,
d\tau\right),\eqno(5.7)$$ with
$K_2=\frac{P(0,0)G_\phi^{n-2}(0)}{\omega_{n-1}(2-n)^2}.$ We can see
that (5.6) coincides with (5.7), the only difference is in constants
$K_1$ and   $K_2.$ Provided initial data satisfy the compatibility
condition
$$G_\phi(0) p_0'(|{ x}|)=-\frac{p(0,0)}{\omega_{n-1}(2-n)^2}\rho_0(|{ x}|)|{ x}|,\eqno(5.8)$$
the density and pressure can be found by formulas (5.4), (5.5).

The analysis of equation (5.6) (or (5.7)) shows that if the constant
$K_1$ (or $K_2$) is positive, as in our case, then the solution
exists for all $t\ge 0$ and $\alpha(t)=O(t^{-1}) $ as
$\,t\to\infty.$

We see that in the second case ($\phi=|{ x}|^{2-n}$) we get a more
wide class of solutions. Indeed, in  Section 5.1 we need to require
the significant decay of solutions as $|{ x}|\to\infty$ to guarantee
the finiteness both of energy and momentum of mass, i.e. we restrict
ourselves by solutions from $K_{GD}$ with $\alpha_v=1.$ When we
construct solutions using the generalized momentum, we do not
require the finiteness of energy and the behavior of smooth pressure
and density, prescribed by convergence of $G_\phi(t)$ and condition
(5.8) is the following: $$p(|{ x}|)=O(|{ x}|^q),\quad q<0, \quad
\rho(|{ x}|)=O(|{ x}|^s),\quad s=q-2<0,\quad |{ x}|\to \infty.$$

\section{Behavior of boundary of a liquid volume}

We mention here very briefly one more application of the generalized
momenta of mass to (GD) equations  (see \cite{RozLV} for details).
Namely, the expansion of boundary of a material volume inside of a
smooth flow of gas can be studied by this method. The last question
is connected with a problem of air pollution:
 one can be interested if under usual meteorological conditions a polluted cloud will attain
some geographical object.

We suppose that initially a point $x_0$ do not belong to ${\mathcal
V}(t).$ Let us set the following question: what conditions we have
to impose on initial data provided they are known only inside of
${\mathcal V}(0),$ to guarantee that the boundary of given material
volume within a smooth flow will attain a given $\varepsilon$ -
neighborhood of point $x_0$? It is clear that for the answer to this
question certain assumptions on the thermodynamic variables in the
whole space have to be done. To formulate these assumptions we need
the following definition.

\begin{definition}
We say that the pressure in the moment $t$ is
\underline{distributed} along the boundary $\partial {\mathcal
V}(t)$ of domain ${\mathcal V}(t), \,x_0\notin {\mathcal V}(t),$
\underline{ regularly with constant $M\ge 0,$ } if
$\Bigl|\int\limits_{\partial{\mathcal V}(t)}\left(\frac{{ x}-x_0}
{|{ x}-x_0|},{\nu}\right) p\,d\partial {\mathcal V}\Bigr|\le M,$
where ${\nu}$ is a unit outer normal.
\end{definition}
If $p(t,x)$ is constant, then the integral in the left hand side of
the latter inequality is equal to zero. Setting $M$ sufficiently
small we assume that the material volume will not meet a zone of
large gradient of pressure.

To prove the following theorem we use the generalized momentum of
mass with $\phi(|{ x}|)=|{ x}-x_0|^q,\, q<0.$

\begin{theorem} {\rm \cite{RozLV}}
Let a finite material volume ${\mathcal V}(t)$ of compressible
liquid
 with $C^1$ - smooth boundary $\partial{\mathcal V}(t)$ do not
contain a point $x_0.$ Suppose that the flow is $C^1$ -- smooth for
all $t\in[0,T],\,T\le\infty,$ and the pressure along the boundary
$\partial{\mathcal V}(t)$ is distributed regularly with a constant
$M$ uniformly in $t$ for $t\in[0,T].$

Let us choose some real numbers  $q<-n-\frac{2}{\gamma-1} $ and
$\varepsilon,$  $\,\,0<\varepsilon<{\rm dist}(\partial{\mathcal
V}(0),x_0).$ Then for all initial data there exists such constant
$\delta\le 0,$ depending on initial data, $\varepsilon,
\,q,\,T,\,M,\,n,\,\gamma$ that if initially $$
\int\limits_{{\mathcal V}(0)}|{ x}-x_0|^{q-2}\,({ v}(0,x),
{x}-x_0)\,\rho_0(x)\,\,dx<\delta, $$ then within a time $t_1$, later
then $T,$ the boundary of given material volume will attain  the
$\varepsilon$ -- neighborhood of point $x_0.$
\end{theorem}




\end{document}